\documentclass{amsart}
\usepackage{amsmath}
\usepackage{amssymb}
\usepackage{graphicx}
\usepackage{latexsym}
\usepackage{color}
\usepackage{amscd}
\usepackage[all]{xy}
\usepackage{enumerate}
\usepackage{hyperref}
\usepackage{subfigure}
\usepackage{soul}
\newtheorem{thm}{Theorem}

\newtheorem{theorem}[thm]{Theorem}

\theoremstyle{definition}

\newtheorem*{definition*}{Definition}

\newcommand{\CPb}{\overline{\mathbb{CP}}{}^{2}}
\newcommand{\CP}{{\mathbb{CP}}{}^{2}}

\newcommand{\Z}{\mathbb{Z}}

\newcommand{\M}{\mathrm{Mod}}

\def \x {\times}

\begin{document}

\title[On four-manifolds without $1$-- and $3$--handles ]
{On four-manifolds without $1$-- and $3$--handles }

\author[R. \.{I}. Baykur]{R. \.{I}nan\c{c} Baykur}
\address{Department of Mathematics and Statistics, University of Massachusetts, Amherst, MA 01003-9305, USA}
\email{inanc.baykur@umass.edu}

\begin{abstract}
We note that infinitely many irreducible, closed, simply connected \linebreak $4$--manifolds, with prescribed signature and spin type, admit perfect Morse functions, i.e. they can be given handle decompositions without $1$- and $3$-handles. \linebreak In particular, there are many such $4$--manifolds homeomorphic but not diffeomorphic to the standard $4$--manifolds $\#_m (S^2 \x S^2)$ and $\#_n (\CP \# \CPb)$, respectively,  which answers Problem 4.91 on Kirby's 1997 list.
\end{abstract}

\maketitle

\setcounter{secnumdepth}{2}
\setcounter{section}{0}

\section{Introduction}
 
A Morse function on a closed manifold $X$ is  \emph{perfect} if the number of its critical points is equal to the sum of its Betti numbers.
The Problem 4.91 on Kirby's 1997 list \cite{K2} by Shkolnikov asks whether any closed, simply connected $4$--manifold $X$ with signature zero and admitting a perfect Morse function is diffeomorphic to $\#_m (S^2 \x S^2)$ or $\#_n (\CP \# \CPb)$.  We answer this question negatively, with many spin and non-spin counter-examples. 
Our examples are handed by the following more general result:

\begin{theorem} \label{thm:main}
There exist irreducible, closed, simply connected  \mbox{$4$--manifolds} that admit perfect Morse functions, where the $4$--manifolds have prescribed signature $\sigma$ and prescribed spin type when $\sigma$ is divisible by $16$, and arbitrarily large Euler characteristic. 
\end{theorem}

A closed smooth $4$--manifold is said to be \emph{irreducible} if it is not diffeomorphic to a connected sum of two $4$--manifolds, neither one of which is a homotopy $4$--sphere. \linebreak \emph{Reducible} examples satisfying the remaining hypotheses of the theorem are easy to come by. For example, one can take connected sums of copies of $\CP$, $S^2 \x S^2$, and the $K3$ surface, with either orientation; all of these standard $4$--manifolds admit perfect Morse functions, and so do their connected sums. Our irreducible $4$--manifolds with signature zero in the theorem are homeomorphic to some $\#_m (S^2 \x S^2)$ or $\#_n (\CP \# \CPb)$ by Freedman, but are not diffeomorphic to them.

The Kirby Problem 4.91 can be seen as a probe into the long-standing open question of whether every \emph{simply connected} $4$-manifold admits a perfect Morse function,  or equivalently, a handle decomposition without $1$-- and $3$--handles. By the celebrated works of Smale in higher dimensions, and Perelman in dimension three, a closed, simply connected, smooth $n$--manifold is known to admit a handle decomposition without \linebreak $1-$ and $(n-1)$--handles in any dimension $n \neq 4$. While many exotic copies of $\#_m (S^2 \x S^2)$ and $\#_n (\CP \# \CPb)$ have been constructed over the past two decades,  inferring information on their handle decompositions seems out of reach for the majority of examples. (This is largely due to their constructions involving compact complex surfaces with non-negative signatures as key building blocks but without a proper understanding of their smooth topology.) The examples we provide in this note are derived from the relatively recent work in  \cite{baykur-hamada:arbitrarysignatureLF}, which, through Lefschetz fibrations,  allows us to present the desired \mbox{$4$--manifolds} with perfect Morse functions

\section{Proof of the theorem and examples}

We first outline an abstract construction scheme for producing exotic simply connected signature zero $4$--manifolds with perfect Morse functions via Lefschetz fibrations. We then explain how to cook up such examples with particular ingredients. After that, we explain how to prove the more general statement given in Theorem~\ref{thm:main}.

Below we assume that the Lefschetz fibrations have non-empty critical locus and are \emph{relatively minimal}, i.e. none of the vanishing cycles are null-homotopic on the fiber. For more details on Lefschetz fibrations and the conventions we adopt in this note, we refer the reader to \cite{baykur-hamada:arbitrarysignatureLF}. 

\medskip
\noindent \underline{\textit{The recipe}}: Let $F \cong \Sigma_g$ be a regular fiber of a Lefschetz fibration $X \to S^2$. There is an induced handle decomposition on $X \setminus \nu F$,  where one takes a perfect Morse function on $D^2 \x \Sigma_g$,  and for each Lefschetz critical point,  attaches a \mbox{$2$--handle} along the corresponding vanishing cycle. Hereon, let $\{a_i, b_i\}_{i=1}^{g}$ denote the standard geometric basis of $H_1(\Sigma_g)$, using which we can label the $1$--handles of the above handle decomposition.

Let $X \to S^2$ be a genus--$g$ Lefschetz fibration with signature $\sigma(X)=0$. Note that by Endo's signature formula for hyperelliptic fibrations \cite{endo}, $\sigma(X) \geq 0$ implies that $g \geq 3$. The fibration has at least one \emph{non-separating} vanishing cycle by \cite{smith}; call it $c_1$.  Let $\phi_i$ and $\psi_i$ be elements of the mapping class group $\M(\Sigma_g)$ such that $\phi_i(c_1)=a_i$ and $\psi_i(c_1)=b_i$ \, for $i=1, \ldots, 2g$. Let $P=1$ in $
M(\Sigma_g)$ be a monodromy factorization for this Lefschetz fibration, where $P:=\Pi_{i=1}^\ell t_{c_i}$ is a product of $\ell$ positive Dehn twists.

We can describe a new Lefschetz fibration $Y \to S^2$ with monodromy factorization
\[
P^{\phi_1} P^{\psi_1} P^{\phi_2} P^{\psi_2} \cdots  P^{\phi_g} P^{\psi_g}=1 \ \ \text{ in } \M(\Sigma_g) \, ,
\]
where $P^\varphi$ for any $\varphi \in \M(\Sigma_g)$ denotes the conjugated product $\Pi_{i=1}^\ell t_{\varphi(c_i)}$. The above monodromy factorization is Hurwitz equivalent to
\[
t_{a_1} t_{b_1} t_{a_2} t_{b_2} \cdots t_{a_g} t_{b_g} P_0 =1 \ \ \text{ in } \M(\Sigma_g) \, ,
\]
where $P_0$ is a product of $\ell_0:=2g (\ell -1)$ positive Dehn twists.  In turn, we get a Lefschetz fibration $Z \to S^2$ with monodromy factorization
\begin{equation} \label{eq:signzero-monodromy}
(t_{a_1} t_{b_1} t_{a_2} t_{b_2} \cdots t_{a_g} t_{b_g} P_0)^2 =1 \ \ \text{ in } \M(\Sigma_g) \, .
\end{equation}
We claim that $Z$ is a simply connected, irreducible $4$--manifold that admits a perfect Morse function.  

The Lefschetz fibration $Z$ can be equipped with a Gompf-Thurston symplectic form, and because it is a fiber sum of (non-trivial) Lefschetz fibrations, the total space $Z$ is minimal by \cite{usher, baykur:LFminimality}. Since minimal rational or ruled surfaces do not admit non-trivial relatively minimal Lefschetz fibrations, we conclude two things: the minimal symplectic $4$--manifold $Z$ is irreducible \cite{hamilton-kotschick},  and it cannot be diffeomorphic to any $\#_m (S^2 \x S^2)$ or $\#_n (\CP \# \CPb)$.

Moreover, we have a (fiber sum) decomposition $Z = (Y \setminus \nu F) \cup (Y \setminus \nu F)$, where each  $Y \setminus \nu F$ is a Lefschetz fibration over $D^2$ prescribed by the product \, $t_{a_1} t_{b_1} t_{a_2} t_{b_2} \cdots t_{a_g} t_{b_g} P_0$ in $\M(\Sigma_g)$. Here, each copy of $Y \setminus \nu F$ admits a handle decomposition with only one $0$--handle and $\ell_0$ $2$--handles. This follows from taking the standard handle decomposition induced by the Lefschetz fibration $Y \setminus \nu F \to D^2$, and attaching the first $2g$ curves along $\{a_i, b_i\}$; we get a geometrically dual $2$--handle to each $1$--handle, so we can cancel all these handle pairs. After the cancellations, viewing the second copy of $Y \setminus \nu F$ with its handle decomposition turned upside down, we get the desired handle decomposition for $Z$ with no $1$-- or $3$-- handles and $2 \ell_0$ many $2$--handles. So, $Z$ admits a perfect Morse function. As there are no $1$--handles, $Z$ is simply connected. 

Since the initial monodromy factorization $P=1$ has signature zero, so does the product of its conjugates. It follows that $Z$ is a simply connected $4$--manifold with signature $\sigma(Z)=0$.

\medskip
\noindent \underline{\textit{The examples}}: The main ingredients for the above recipe,  namely, (non-trivial) Lefschetz fibrations over $S^2$ with signature zero, were discovered in \cite{baykur-hamada:arbitrarysignatureLF}. The example we handpick for our construction here is the \emph{spin} genus--$9$ Lefschetz fibration $X \to S^2$ \mbox{in \cite{baykur-hamada:arbitrarysignatureLF}.}  In fact, with this particular example, we can improve on the general recipe above: It suffices to use $6$ conjugates of $P$ (instead of $2g=18$ conjugates, that is) to build a monodromy factorization for a \emph{spin} Lefschetz fibration $Y \to S^2$ with vanishing cycles dual to every $\{a_i, b_i\}$. (This is the exotic $\#_{127} (S^2 \x S^2)$ constructed in \cite[Proof of Corollary~C, p.46--47]{baykur-hamada:arbitrarysignatureLF}.)  
A fiber sum of two copies of $Y \to S^2$ gives us the desired $Z$, homeomorphic but not diffeomorphic to $\#_{271} (S^2 \x S^2)$, and admits a perfect Morse function as we had explained. 

Non-spin examples are certainly easier to build, but for the brevity of this note, let us derive them  from the same genus--$9$ fibration. For this, in the final step of our construction above, we take a twisted fiber sum of two copies of $Y \to S^2$ with a mapping class $\varphi \in \M(\Sigma_9)$ mapping $a_9$ curve to $b_9$; take e.g. $\varphi= t_{b_9} t_{a_9}$. This results in a Lefschetz fibration $Z \to S^2$ with the same properties as earlier, except it is not spin. To verify this last claim, note that the monodromy factorization of $Z$ is Hurwitz equivalent to
\[
t_{b_1} t_{b_1'} t_{b_9} P_1 =1 \ \  \text{ in } \M(\Sigma_9) \, ,
\]
where the disjoint curves $b_1, b'_1, b_9$ (which are the curves $\beta_1, \beta'_1, \beta_9$ in \cite[Fig.29]{baykur-hamada:arbitrarysignatureLF}) cobound a genus--$0$ subsurface on $\Sigma_9$, and $P_1$ is the product of the remaining Dehn twists. Perturbing the fibration on $Z$, we get another Lefschetz fibration where these three nodes are clustered on a reducible fiber containing a genus--$0$ component with self-intersection $-3$. So we have a non-spin $Z$ homeomorphic but not diffeomorphic to $\#_{271} (\CP \# \CPb)$.

It should be clear at this point that by starting with  signature zero Lefschetz fibrations of other genera and/or taking further fiber sums in the construction of $Z$, we can generate infinitely many examples with arbitrarily large Euler characteristic. 

\medskip
\noindent \underline{\textit{Any signature}}: To prove the more general statement given in the theorem,  pick any genus--$g$ Lefschetz fibration $X' \to S^2$ with prescribed signature and spin type. Make sure also to have a genus--$g$ Lefschetz fibration $X \to S^2$ with $\sigma(X)=0$ and with monodromy that has a \emph{matching spin type}. Such $X\to S^2$ exists for $g=8d+1$, for any $d \in \Z^+$ by \cite[Theorem~A]{baykur-hamada:arbitrarysignatureLF}.

As before, we build a signature zero Lefschetz fibration $Z \to S^2$ with monodromy factorization~\eqref{eq:signzero-monodromy}.  If $P'=1$ in $\M(\Sigma_g)$ is the monodromy factorization of $X' \to S^2$, then we can obtain a new fibration $Z' \to S^2$ with monodromy factorization 
\[
(t_{a_1} t_{b_1} t_{a_2} t_{b_2} \cdots t_{a_g} t_{b_g} P_0)^2 \, P' 
=1 \ \ \text{ in } \M(\Sigma_g) \, ,
\]
which is Hurwitz equivalent to the monodromy factorization
\begin{equation} \label{eq:arbsignaturemonodromy}
(t_{a_1} t_{b_1} t_{a_2} t_{b_2} \cdots t_{a_g} t_{b_g} P_0) \, P' \,
(t_{a_1} t_{b_1} t_{a_2} t_{b_2} \cdots t_{a_g} t_{b_g} P_0) 
=1 \ \ \text{ in } \M(\Sigma_g) \, .
\end{equation}

A technical side note: In the spin case, we instead express the monodromy factorization in a spin mapping class group $\M(\Sigma_g; s)$ for some spin structure $s \in \mathrm{Spin}(\Sigma_g)$. For $g=8d+1$, there is an $X \to S^2$ with monodromy factorization $P$ in the desired $\M(\Sigma_g; s)$.  This is because these examples of \cite{baykur-hamada:arbitrarysignatureLF} admit spin monodromy factorizations with either Arf invariant; conjugating the one in the right orbit, we can get $P$ to be in the desired spin mapping class group $\M(\Sigma_g; s)$.

It follows that $Z'$ has the same signature and type as $X'$. Here $Z'\to S^2$ is the fiber sum of three Lefschetz fibrations: two copies of $Y \to S^2$ (yielding $Z$) and $X' \to S^2$. Express this fiber sum as $(Y \#_F X') \#_F Y$. By the same arguments as before, $Z'$ is irreducible, and we deduce from the decomposition $Z' = ((Y \#_F X') \setminus \nu F) \cup (Y  \setminus \nu F)$ that it admits a perfect Morse function. For examples with arbitrarily large Euler characteristic, we take fiber sums of $Z' \to S^2$ with further copies of $Y \to S^2$. 
\qed

\bigskip
In this note, we aimed to keep our arguments as simple as possible rather than the underlying topology of the examples. It is plausible---but not at all clear---that the explicit Kirby diagrams one can get for the smallest known examples of exotic $\#_{11} (S^2 \x S^2)$ and $\#_{9} (\CP \# \CPb)$ more recently given in \cite{baykur-hamada:exoticsignaturezero} can also be further simplified to get handle diagrams without $1$-- and $3$--handles.

\bigskip
\noindent \textit{Acknowledgements}: 
The author would like to thank Tye Lidman for helpful comments on a draft of this note. This work was supported by the NSF grant DMS-2005327.

\bigskip
\bibliography{references}

\begin{thebibliography}{1}

\bibitem{baykur:LFminimality}
R.~I. Baykur.
\newblock Minimality and fiber sum decompositions of {L}efschetz fibrations.
\newblock {\em Proc. Amer. Math. Soc.}, 144(5):2275--2284, 2016.
\newblock \href {https://doi.org/10.1090/proc/12835}
  {\path{doi:10.1090/proc/12835}}.

\bibitem{baykur-hamada:arbitrarysignatureLF}
R.~I. Baykur and N.~Hamada.
\newblock Lefschetz fibrations with arbitrary signature, 2020.
\newblock to appear: J. Eur. Math. Soc.
\newblock \href {http://arxiv.org/abs/2010.11916} {\path{arXiv:2010.11916}}.

\bibitem{baykur-hamada:exoticsignaturezero}
R.~I. Baykur and N.~Hamada.
\newblock Exotic 4-manifolds with signature zero, 2023.
\newblock \href {http://arxiv.org/abs/2305.10908} {\path{arXiv:2305.10908}}.

\bibitem{endo}
H.~Endo.
\newblock Meyer's signature cocycle and hyperelliptic fibrations.
\newblock {\em Math. Ann.}, 316(2):237--257, 2000.
\newblock \href {https://doi.org/10.1007/s002080050012}
  {\path{doi:10.1007/s002080050012}}.

\bibitem{hamilton-kotschick}
M.~J.~D. Hamilton and D.~Kotschick.
\newblock Minimality and irreducibility of symplectic four-manifolds.
\newblock {\em Int. Math. Res. Not.}, pages Art. ID 35032, 13, 2006.
\newblock \href {https://doi.org/10.1155/IMRN/2006/35032}
  {\path{doi:10.1155/IMRN/2006/35032}}.

\bibitem{K2}
R.~Kirby.
\newblock Problems in low--dimensional topology.
\newblock In W.~Kazez, editor, {\em Geometric Topology}. American Math.\
  Soc./International Press, Providence, 1997.

\bibitem{smith}
I.~Smith.
\newblock Lefschetz fibrations and the {H}odge bundle.
\newblock {\em Geom. Topol.}, 3:211--233, 1999.
\newblock \href {https://doi.org/10.2140/gt.1999.3.211}
  {\path{doi:10.2140/gt.1999.3.211}}.

\bibitem{usher}
M.~Usher.
\newblock Minimality and symplectic sums.
\newblock {\em Int. Math. Res. Not.}, pages Art. ID 49857, 17, 2006.
\newblock \href {https://doi.org/10.1155/IMRN/2006/49857}
  {\path{doi:10.1155/IMRN/2006/49857}}.

\end{thebibliography}
\bibliographystyle{abbrvurl}

\end{document}